\documentclass[11pt,leqno]{amsart}%
\usepackage{amsmath}
\usepackage{amsfonts}
\usepackage{amssymb}
\usepackage{graphicx}
\usepackage{dsfont}
\usepackage{xcolor}

\usepackage{stix}
\usepackage{hyperref}
\usepackage{bm}
\usepackage[margin=1.90cm]{geometry}

\providecommand{\U}[1]{\protect\rule{.1in}{.1in}}
\newtheorem{theorem}{Theorem}[section]
\newtheorem*{acknowledgement*}{Acknowledgement}

\newtheorem{definition}[theorem]{Definition}

\newtheorem{problem}[theorem]{Problem}

\newtheorem{proposition}[theorem]{Proposition}
\newtheorem{remark}[theorem]{Remark}

\newcommand{\R}{\mathbb{R}}

\newcommand{\p}{\partial}

\newcommand{\Hess}{\mathrm{Hess}}

\title[Compactly supported Sobolev approximations]{Sobolev functions without compactly supported approximations}
\author {Giona Veronelli}
\address[Giona Veronelli]{Dipartimento di Matematica e Applicazioni, Universit\`a di Milano Bicocca, via R. Cozzi 53, I-20126 Milano, Italy}
\email{giona.veronelli@unimib.it}

\begin{document}
\begin{abstract}
A basilar property and a useful tool in the theory of Sobolev spaces is the density of smooth compactly supported functions in the space $W^{k,p}(\R^n)$ (i.e. the functions with weak derivatives of orders $0$ to $k$ in $L^p$). On Riemannian manifolds, it is well known that the same property remains valid under suitable geometric assumptions. However, on a complete non-compact manifold it can fail to be true in general, as we prove in this paper. This settles an open problem raised for instance by E. Hebey [\textit{Nonlinear analysis on manifolds: Sobolev spaces and inequalities}, Courant Lecture Notes in Mathematics, vol. 5, 1999, pp. 48-49].
\end{abstract}

\date{\today}

\maketitle

\section{Introduction}\label{Intro}

Let $(M^{n}, g)$ be a smooth, complete (possibly non-compact) $n$-dimensional Riemannian manifold without boundary, $k$ a natural integer, and $p\geq 1$. Generalizing some of the usual equivalent definitions introduced for $M=\R^n$, one can consider different notions of Sobolev spaces on $M$. Unlike the Eucildean setting, on a Riemannian manifold these definitions a priori may differ from each other. Accordingly, in the last decades a lot of effort has been made to understand the relations between the various notions of Sobolev space on a manifold. This paper aims at giving a further contribution in this direction.

A first definition of Sobolev space involves weak (covariant) derivatives. \begin{definition}\label{def_1}
	The Sobolev space\footnote{From now on, we will omit the dependence on the Riemannian metric in the names of Sobolev spaces. For instance $W^{k,p}(M)$ stands for $W^{k,p}(M,g)$.} $W^{k,p}(M)$ is the space of functions on $M$ all of whose (distributional) covariant derivatives of order $j$ are tensor fields with finite $L^{p}$-norm, for $0\leq j\leq k$. This turns out to be a Banach space, once endowed with the natural norm
\[
 \left\|u\right\|_{W^{k,p}}\dot =\sum_{j=0}^{k}\left(\int_{M}|\nabla^{j}u|^{p}\,\,d\!\mu_g\right)^{\frac{1}{p}}.
\]
\end{definition}
By a generalised Meyers-Serrin-type theorem (see e.g. \cite{GGP-AccFinn}) it is known that on any complete Riemannian manifold $M$, the space $W^{k,p}(M)$ is the closure of 
\[
\mathcal{C}_{k}^{p}(M,g)\doteq C^{\infty}(M)\cap W^{k,p}(M,g)=\left\{u\in C^{\infty}(M): \int_{M}|\nabla^{j}u|^{p}\,d\!\mu_g<+\infty\,\,\forall\,j=0,\ldots,k\right\}
\]
with respect to the norm $\left\|\cdot\right\|_{W^{k,p}(M)}$.
In the applications, it is often very useful to prove a result in the class of compactly supported functions, and hence to extend its validity to the whole Sobolev space by a density argument. This leads to a second definition of Sobolev spaces.
\begin{definition}\label{def_3}
	The space $W_{0}^{k,p}(M)\subseteq W^{k,p}(M)$ is the closure of the space of smooth compactly supported functions $C_{c}^{\infty}(M)$ with respect to the norm $\left\|\cdot\right\|_{W^{k,p}}$.
	\end{definition}

The main result of this article will show that, in all generality, the inclusion $W_{0}^{k,p}(M)\subseteq W^{k,p}(M)$ can be a proper inclusion (see Theorem \ref{th_main} below). In order to contextualize the problem, let us introduce some other spaces.
A third natural definition of Sobolev space on $M$, which in the special case $M=\R^n$ is modelled on operator theory or Fourier transform, is the following.
\begin{definition}\label{def_4}
	 For even $k=2m\in\mathbb N$, we define the space
\[
H^{k,p}(M)\dot = \left\{ u\in L^p(M)\,:\, (1-\Delta)^m u\in L^p \right\},
\]
where $\Delta$ is the Laplace-Beltrami operator of $M$ (with negative spectrum) and the derivatives have to be interpreted in the distributional sense. The natural associated norm is 
\[ \left\|u\right\|_{H^{k,p}}\dot=\|(1-\Delta)^m u\|_{L^p}.
\]
\end{definition}
Inspired to Definitions \ref{def_1} and \ref{def_3}, it is natural to consider also the following spaces.
\begin{definition}\label{def_5}
	For even $k\geq 2$, we define the spaces
	\[
	\widetilde H^{k,p}(M)\dot =  H^{k,p}(M)\cap  W^{k-1,p}(M).
	\]
\end{definition}
\begin{definition}
	For even $k\geq 2$,	the space $H_{0}^{k,p}(M)\subseteq H^{k,p}(M)$ is the closure of $C_{c}^{\infty}(M)$ with respect to the norm $\left\|\cdot\right\|_{H^{k,p}}$.
\end{definition}
 
During the last decades, several efforts have been made in order to understand which are the relations and inclusions between the spaces above. It is by now well known that all these definitions coincide when $M$ is $\R$, or is a closed (i.e. compact without boundary) manifold. The situation is much more complicated when $M$ is complete non-compact. 
Clearly, by their very definitions, one has the inclusions
\[
W_0^{k,p}(M,g)\subseteq W^{k,p}(M,g)\subseteq \widetilde H^{k,p}(M,g)\subseteq  H^{k,p}(M,g),
\]
for every $p\in[1,\infty)$ and any integer $k$ (assuming $k$ even in the second and third inclusion).

Consider the operator $\mathcal J=(1-\Delta)^{m}$ defined on $C^\infty_c(M)$. We notice here that $H^{2m,2}_0(M)$ and $H^{2m,2}(M)$ are respectively the domain of the closure of $\mathcal J$ and the domain of the adjoint of $\mathcal J$ in $L^2(M)$.
 As a consequence of Chernoff's theorem, i.e. the essential self-adjointness of $(1-\Delta)^m$ on $C^\infty_c(M)$, \cite{chernoff}, one gets that 
\[
H_0^{k,2}(M,g)=  H^{k,2}(M,g),
\]
on any complete manifold for every even integer $k$.
At least for $k=2$ one has also that 
\[
H_0^{2,p}(M,g)= H^{2,p}(M,g)\] 
for every $p\in(1,\infty)$, as proved by O. Milatovic in an appendix of \cite{GuneysuPigola_AMPA}. Moreover,
using \cite[Corollary 3.5]{stric} (for $p=2$) or \cite[Theorem 4.1]{CD} one gets also
\[
\tilde H^{2,p}(M,g)= H_0^{2,p}(M,g)=H^{2,p}(M,g)\]
for every $1<p\leq 2$.
Conversely, it turns out that 
\begin{equation}\label{dod}
W^{2,p}(M)\subsetneq H^{2,p}(M)
\end{equation}
 with proper inclusion. This fact is known at least since \cite[p. 67]{dodziuk}, where J. Dodziuk proposed a counterexample in the case of Sobolev spaces of $1$-forms on a complete surface for $p=2$. His example exploits the Weitzenb\"ock identity applied to harmonic 1-forms. It seems thus hardly generalizable to Sobolev spaces of functions, for which we are not aware of any counterexample in the literature. We will give such a counterexample in Section \ref{app}.\\
  
In this paper we address in particular the following problem.

\begin{problem}\label{Pb}
	Given an arbitrary complete (non-compact) Riemannian manifold $(M,g)$, is it true that
	\begin{equation}\label{eq_dens}
	W_{0}^{k,p}(M)=W^{k,p}(M)
	\end{equation} 
	for all integer $k\geq 0$ and $p\in[1,\infty)$?
\end{problem}

This problem has been intensively studied, and several partial results are known. For instance, it is a standard fact that $W_{0}^{0,p}(M)=W^{0,p}(M)=L^p(M)$, and with a little effort one can also prove that $W_{0}^{1,p}(M)=W^{1,p}(M)$ for all $k$ and $p\in[1,\infty)$ on any complete manifold; \cite{aubin-bull}. Concerning the non-trivial case $k\geq 2$, to the best of our knowledge the most general and up-to-date result is the following theorem from \cite{IRV-HO}, which generalizes previously known achievements treating the case of constant bounds on the curvature, or specific to the 2nd order $k=2$. A non-exhaustive list of older results in this direction includes contributions by T. Aubin \cite{aubin-bull}, E. Hebey \cite{Hebey,HebeyCourant}, L. Bandara \cite{Bandara}, B. G\"uneysu and S. Pigola \cite{Guneysu-Book},  \cite{GuneysuPigola}, and by the author in collaboration with D. Impera and M. Rimoldi, \cite{IRV-HessCutOff}.

\begin{theorem}[Theorem 1.5 in \cite{IRV-HO}]\label{th_IRV}
	Let $(M,g)$ be a complete Riemannian manifold and $r(x)\dot=\operatorname{dist}_g(x,o)$ the Riemannian distance function from a fixed reference point $o\in M$. 
Define \[\lambda(r(x))\dot=r(x)\prod_{j=1}^{K}\ln^{[j]}(r(x)),\] where $\ln^{[j]}$ stands for the $j$-th iterated logarithm (e.g. $\ln^{[2]}(t)=\ln\ln t$, etc.) and $K$ is some positive integer.

We have that  $W^{k,p}(M)=W^{k,p}_0(M)$ for all $p\in[1,+\infty)$ if 	\[
	\ |\nabla^j\mathrm{Ric}|(x)\leq \lambda(r(x))
	^{(2+j)/(k-1)},\qquad 0\leq j \leq k-2,\] 
	and either
	\begin{itemize}
		\item[(a)] for some $i_0>0$,
		$\mathrm{inj}(x)\geq i_0\lambda(r(x))^{-1/(k-1)}>0$, or
		\item[(b)] for some $D>0$,
		$|\mathrm{Riem}|(x)\leq D^2\lambda(r(x))^{2/(k-1)}.$
	\end{itemize}
Moreover $W^{2,2}(M)=W^{2,2}_0(M)$ if 
\[
\ \mathrm{Ric}(x)\geq -\lambda(r(x))
^{2},\]
and $W^{k,2}(M)=W^{k,2}_0(M)$ for $k>2$ if 	\[
\ |\nabla^j\mathrm{Riem}|(x)\leq \lambda(r(x))^{(2+j)/(k-1)},\qquad 0\leq j \leq k-3.\]
	
\end{theorem} 

Accordingly, on any Riemannian manifold whose metric is (slightly) controlled, compactly supported functions are dense in $W^{k,p}(M,g)$ and Problem \ref{Pb} has a positive answer. It seems thus very natural to ask if this fact remains true in full generality (i.e on any complete Riemannian manifold), or 
if the inclusion $W_0^{k,p}(M,g)\subseteq W^{k,p}(M,g)$ turns out to be a proper inclusion for some very wild geometries. 
Quite surprisingly, Problem \ref{Pb} apparently remained open in general so far; see for instance \cite[pp. 48-49]{HebeyCourant}. 
This paper aims at filling this gap.

\begin{theorem}\label{th_main}
	For any integer $n\geq 2$, there exists a complete $n$-dimensional Riemannian manifold $(M.g)$ such that for any integer $k\geq 2$ and any real $p\geq 2$, the space of smooth compactly-supported functions $C^\infty_c(M)$ is not dense in $W^{k,p}$ with respect to the norm $\|\cdot\|_{W^{k,p}}$. 
\end{theorem}

Problem \ref{Pb} is quite subtle because of the following phenomenon. Let $k=p=2$. Combining Chernoff's Theorem, \cite{chernoff}, with an $L^2$-gradient estimate, \cite[Corollary 3.5]{stric}, one gets that, on any complete Riemannian manifold $M$, the class of functions $C^\infty_c(M)$ is dense in $\tilde H^{2,2}(M)$ (hence in $W^{2,2}(M)$) with respect to the norm $\left\|\cdot\right\|_{L^{2}}+\left\|\nabla \cdot\right\|_{L^2}+\left\|\Delta \cdot\right\|_{L^{2}}$.
Accordingly, in order to prove Theorem \ref{th_main}, one has to find a manifold $(M,g)$ and a function $F\in W^{2,2}(M)$ such that
\[
 \|\eta-F\|_{W^{2,2}(M)}=\|\eta-F\|_{W^{1,2}(M)}+\|\nabla^2(\eta-F)\|_{L^{2}(M)}>\epsilon_0\]
for every $\eta\in C^\infty_c(M)$ and for some $\epsilon_0>0$ independent of $\eta$, whereas $\|\eta-F\|_{W^{1,2}(M)}+\|\Delta(\eta-F)\|_{L^{2}(M)}$ can be made arbitrarily small. 

To achieve this purpose, we proceed as follows. 
For simplicity, choose $\R	\times \mathbb S^{n-1}$ as the underlying space for $(M,g)$, and consider Riemannian metrics which are rotationally symmetric, i.e. $g=dt^2+j^2(t)\sigma$, with $j\in C^\infty(\R,(0,+\infty))$ and $\sigma$ the standard round metric on $\mathbb S^{n-1}$ of constant sectional curvature $1$. In particular we can make $F$ assume different constant values (says $0$ and $1$) when $t$ is respectively small enough or large enough. If $j$ decays fast enough so that $(M,g)$ has  finite volume, then $F$ is in $W^{k,2}(M)$ since it is eventually constant. Now we have to choose the (rapidly decaying) warping function $j$ so that $F$ can not be approached by compactly supported functions in the $W^{2,2}$-topology. The key ingredient is a Hardy-type inequality satisfied by the gradient of $W^{2,2}(M)$ functions, that is
\begin{equation}\label{hardy}
\|\nabla \phi\|_{W^{1,2}(M)} \geq \int_M |\nabla \hat \phi|^2(x)\,\omega(x) \,d\!\mu_g(x),   
\end{equation}
where $\hat \phi$ is the radial symmetrization of $F$ which depends only on the $\R$ component of $M=\R\times \mathbb S^{n-1}$, and $\omega$ is a rotationally symmetric weight which depends on the metric $g$ of $M$ and on its derivatives, namely \[\omega(x)=j(t)^{n-1} + (n-1)j(t)^{n-3}(j'(t))^2.\] 
By our choice for $F$, for every $\eta\in C^\infty_c(M)$, $\phi\dot=F-\eta$ takes different constant values on each end of $M$. In particular, $\nabla \hat \phi$ can not be arbitrarily small. More precisely we can prove a uniform (i.e. independent of $\eta\in C^\infty_c(M)$) lower bound $\epsilon_0>0$ on the RHS of \eqref{hardy}, provided that $|j'|$ grows fast enough in a set of large measure.
Accordingly, it is enough to let $j$ decay uniformly at infinity (so that $M$ has finite volume), but with sufficiently frequent oscillations. \\	

It is worth to observe that the obstruction to get $W_0^{k.p}(M)=W^{k.p}(M)$ is not of topological nature. Even if, in order to simplify technicalities, we construct our counterexample on the non-contractible cylinder $\R\times\mathbb S^{n-1}$, a similar argument can be adapted to get a complete metric on $\R^{n}$. Indeed, it is enough to consider $(M,g)=(\R\times\R^{n-1},dt^2+j^2(t)\bar g)$, with $\bar g$ any complete metric of finite volume on $\R^{n-1}$.

The density of compactly supported functions in $W^{2,k}(M)$ is intimately related to the validity of a global Calder\'on-Zygmund (C-Z in short) inequality on $M$, i.e. $\|\Hess \varphi\|_{L^p}\leq \|\varphi\|_{H^{2,p}}$ for any $\varphi\in C^\infty_c(M)$. Indeed, as we will explain in Remark \ref{rmk_CZ}, examples of
manifolds such that $W^{2,k}_0(M)\neq W^{2,k}(M)$
are automatically also examples of manifolds on which a C-Z inequality can not hold (see \cite{GuneysuPigola,Li} for an introduction to the C-Z  inequalities and previous counterexamples). However, the converse is not true in general, so that giving a negative answer to Problem \ref{Pb} seems, in a sense, harder than disproving C-Z inequalities on some $M$. In this regard, note that the assumptions guaranteeing the validity of a C-Z inequality require some uniform bound on the curvature, and are thus stronger  than the assumptions of Theorem \ref{th_IRV}. \\

In the next section we will prove Theorem \ref{th_main}. Then, in the short final section we will prove the proper inclusion \eqref{dod}.

\section{Proof of Theorem \ref{th_main}}

Let $(M,g):=(\R\times\mathbb S^{n-1},dt^2+j^2(t)\sigma)$, with $j\in C^\infty(\R,(0,+\infty))$ and $\sigma$ the standard round metric on $\mathbb S^{n-1}$ of constant sectional curvature $1$. Note that $(M,g)$ is complete as long as $j>0$.

Before choosing the function $j$ which gives the counterexemple, let us observe some general fact.
Let $(\theta_2,\dots,\theta_n)$ be a local coordinate system of $\mathbb S^{n-1}$. Denote $\partial_i:=\partial_{\theta_i}$ and suppose that
$(\partial_2,\dots,\partial_n)$ is
orthonormal with respect to $\sigma$ at a point $P\in \mathbb S^{n-1}$. Then $\{\partial_t\}\cup \{j(t)^{-1}\partial_i\}_{i=2}^n$ is a local coordinate frame of $M$ orthonormal with respect to the warped metric $g$ at the point $(t,P)\in M$. We have
\begin{align*}
&\Hess_g\,F|_{(t,P)}(\partial_t,\partial_t)=\partial^2_{tt}F(t,P)\\
&\Hess_g\,F|_{(t,P)}(\partial_t,\partial_i)=\partial^2_{ti}F(t,P)-j(t)^{-1}j'(t)\p_iF(t,P)\\
&\Hess_g\,F|_{(t,P)}(\partial_i,\partial_k)=\Hess_{\sigma}\,(F(t,\cdot))|_{P}(\partial_i,\partial_k)-j(t)j'(t)\delta_{ik}\p_tF(t,P),
\end{align*}
for all $i,k=2,\dots,n$. Hence
\begin{align*}
|\Hess_g\,F|_{(t,P)}|_g^2 &= (\Hess_g\,F|_{(t,P)}(\partial_t,\partial_t))^2 + \sum_{i=2}^n (\Hess_g\,F|_{(t,P)}(\partial_t,\partial_i))^2 j(t)^{-2} + \sum_{i,k=2}^n (\Hess_g\,F|_{(t,P)}(\partial_i,\partial_k))^2 j(t)^{-4}\\
&= (\partial^2_{tt}F(t,P))^2 + \frac{|\nabla_\sigma \partial_tF|_\sigma^2(t,P)}{j^2(t)} + \frac{(j'(t))^2}{j(t)^4}|\nabla_\sigma F|_\sigma^2(t,P) -2 \frac{j'(t)}{j(t)^3} \sigma(\nabla_\sigma \p_t F(t,P), \nabla_\sigma  F(t,P)) \\
&+ \frac{|\Hess_\sigma\,F|_{(t,P)}|_\sigma^2}{j^4(t)}+ (n-1)\frac{(j'(t))^2}{j(t)^2}|\p_t F(t,P)|^2 - 2 \frac{j'(t)}{j(t)^3} \Delta_\sigma F(t,P)\p_tF(t,P).
\end{align*}
Integrating on $\mathbb S^{n-1}$ we thus have that for any $t\in \mathbb R$
\begin{align*}
\int_{\mathbb{S}^{n-1}}|\Hess_g\,F|_{(t,P)}|_g^2\,\,d\!\mu_\sigma(P) &\geq  \int_{\mathbb{S}^{n-1}} (\partial^2_{tt}F(t,P))^2 \,\,d\!\mu_\sigma(P)
 + (n-1) \int_{\mathbb{S}^{n-1}} \frac{(j'(t))^2}{j(t)^2}|\p_t F(t,P)|^2 \,\,d\!\mu_\sigma(P)\\
& -2 \frac{j'(t)}{j(t)^3}\int_{\mathbb{S}^{n-1}} \left[ \sigma(\nabla_\sigma \p_t F(t,P), \nabla_\sigma  F(t,P)) + \Delta_\sigma F(t,P)\p_tF(t,P)\right] \,\,d\!\mu_\sigma(P)\\
&=  \int_{\mathbb{S}^{n-1}} (\partial^2_{tt}F(t,P))^2 \,\,d\!\mu_\sigma(P)
+ (n-1) \int_{\mathbb{S}^{n-1}} \frac{(j'(t))^2}{j(t)^2}|\p_t F(t,P)|^2 \,\,d\!\mu_\sigma(P)
\end{align*}
where we have used the Stokes' theorem in the last equality. Also,
\[
|\nabla_g F|_g^2(t,P) = (\p_t F (t,P))^2 + \sum_{i,k=2}^nj^2(t)\sigma (\p_i F (t,P),\p_k F (t,P))  = (\p_t F (t,P))^2 + j^2(t)|\nabla_\sigma F(t,\cdot)|_\sigma^2(P)
\]
so that
\[
\int_{\mathbb{S}^{n-1}}|\nabla_g F|_g^2(t,P)\,\,d\!\mu_\sigma(P) \geq \int_{\mathbb{S}^{n-1}}(\p_t F (t,P))^2\,\,d\!\mu_\sigma(P) 
\]
Define now 
\begin{equation}\label{hat}
\hat F(t):=  \intbar_{\mathbb{S}^{n-1}} F(t,P)\,\,d\!\mu_\sigma(P):= \frac{1}{\omega_n}\int_{\mathbb{S}^{n-1}} F(t,P)\,\,d\!\mu_\sigma(P),
\end{equation}
where $\omega_n:=\mathrm{vol}_\sigma(\mathbb{S}^{n-1})$.
By Jensen inequality 
\[
\intbar_{\mathbb{S}^{n-1}} \hat F^2(t)\,\,d\!\mu_\sigma(P)=\hat F(t)^2\leq \intbar_{\mathbb{S}^{n-1}} F^2(t,P)\,\,d\!\mu_\sigma(P).
\]
Similarly
\[
\intbar_{\mathbb{S}^{n-1}} (\p_t \hat F(t))^2\,\,d\!\mu_\sigma(P)= (\p_t \hat F(t))^2= \left(\intbar_{\mathbb{S}^{n-1}} \p_t  F(t,P)\,\,d\!\mu_\sigma(P) \right)^2 
\leq \intbar_{\mathbb{S}^{n-1}} (\p_tF(t,P))^2\,\,d\!\mu_\sigma(P),
\]
and 
\[
\intbar_{\mathbb{S}^{n-1}} (\p^2_{tt}\hat F(t))^2\,\,d\!\mu_\sigma(P)=(\p^2_{tt} \hat F(t))^2= \left(\intbar_{\mathbb{S}^{n-1}} \p^2_{tt}  F(t,P)\,\,d\!\mu_\sigma(P) \right)^2 
\leq \intbar_{\mathbb{S}^{n-1}} (\p^2_{tt}F(t,P))^2\,\,d\!\mu_\sigma(P).
\]
Hence
\begin{align*}
\|F\|^2_{W^{2,2}(M)} &= \int_{-\infty}^{+\infty} j(t)^{n-1}\int_{\mathbb{S}^{n-1}}\left\{|\Hess_g\,F|_{(t,P)}|_g^2 +|\nabla_g F(t,P)|_g^2 +(F(t,P))^2   \right\}\,\,d\!\mu_\sigma(P)\,dt\\
&\geq  \int_{-\infty}^{+\infty}j(t)^{n-1}\int_{\mathbb{S}^{n-1}}\left\{ (\partial^2_{tt}F(t,P))^2 
+\left[1+ (n-1)  \frac{(j'(t))^2}{j(t)^2}\right](\p_t F(t,P))^2+F (t,P)^2\right\} \,\,d\!\mu_\sigma(P)\,dt\\
&\geq \int_{-\infty}^{+\infty}j(t)^{n-1}\int_{\mathbb{S}^{n-1}}\left\{ (\partial^2_{tt}\hat F(t))^2 
+\left[1+ (n-1)  \frac{(j'(t))^2}{j(t)^2}\right](\p_t \hat F(t))^2+\hat F (t)^2\right\} \,\,d\!\mu_\sigma(P)\,dt\\
&= \|\hat F\|^2_{W^{2,2}(M)} 
\end{align*}
where, with a slight abuse of notation, we interpret  $\hat F(t)$ also as the rotationally symmetric function on $M$ defined by $(t,P)\mapsto \hat F (t)$. We have thus in particular that 
\begin{align}\label{rearra}
\|F\|^2_{W^{2,2}(M)} 
&\geq \int_{-\infty}^{+\infty}j(t)^{n-1}\int_{\mathbb{S}^{n-1}} 
\left[1+ (n-1)  \frac{(j'(t))^2}{j(t)^2}\right](\p_t \hat F(t))^2 \,\,d\!\mu_\sigma(P)\,dt \\
&\geq \omega_n\int_{-\infty}^{+\infty}
\left[j(t)^{n-1}+ (n-1) j(t)^{n-3}(j'(t))^2\right](\p_t \hat F(t))^2 \,dt \nonumber
\end{align}

Now, let $f\in C^\infty (M)$ a smooth function which satisfy 
\[
\begin{cases}
f(t,P)=1,&\text{on } [1,\infty]\times \mathbb{S}^{n-1}\\
f(t,P)=0,&\text{on } [-\infty,-1]\times \mathbb{S}^{n-1}
\end{cases}
\]
In particular 
\[
\nabla^k_g f=0,\ \text{on}\ ([-\infty,-1]\cup [1,\infty])\times \mathbb{S}^{n-1}
\]
for every $k\ge 1$ (including  $\mathrm{Hess}_g\,f=\nabla^2_gf$), so that $\nabla^k_g f\in  L^p(M)$ trivially. We will choose the warping function $j(t)$ so that 
\begin{equation}\label{jvol}
j^{n-1}\in L^1(\R).
\end{equation}This implies that  $\mathrm{vol}_g(M)<\infty$, and thus $f\in L^p(M)$ and $f\in W^{k,p}(M)$ for any $k\geq 2$ and $p\geq 1$. 

Now, suppose that there exists a family $\{f_k\}_{k=1}^\infty\subset  C^{\infty}_c(M)$ of smooth compactly supported function such that $f_k\to f$ in the $W^{k,p}(M)$-norm as $k\to \infty$, for some $p\geq 2$. In particular one has that 
\begin{equation}\label{approx0}
\phi_k:=f-f_k\stackrel{k\to\infty}{\longrightarrow}\ 0\quad\text{in }W^{2,p}(M).
\end{equation}
Since $p\geq 2$ and $M$ has finite volume, by H\"older inequality one deduces that
\begin{equation}\label{approx}
\phi_k:=f-f_k\stackrel{k\to\infty}{\longrightarrow}\ 0\quad\text{in }W^{2,2}(M).
\end{equation}

Since the $f_k$'s are compactly supported on $M$, there exist constants $t_k>0$ depending on $k$ such that  $\phi_k(t,p)= 1$ on $[t_k,+\infty)\times \mathbb S^{n-1}$ and  $\phi_k(t,p)= 0$ on $(-\infty,-t_k]\times\mathbb{S}^{n-1}$.

To get a contradiction with \eqref{approx}, we are going to prove that, for a suitable choice of the warping function $j(t)$, there exists a strictly positive constant $C_j$ depending on $j(t)$, but independent of $k$, such that $\|\phi_k\|_{W^{2,2}(M)}\geq C_j>0$. To this end, note first that $\hat \phi_k(t)= 1$ on $[t_k,+\infty)$ and $\hat\phi_k(t)= 0$ on $(-\infty,-t_k]$, where the function $\hat\phi_k$ with a superscript $\hat{\ }$ is defined as in \eqref{hat}. In particular 
\[
\int_{-\infty}^{\infty} \p_t\hat\phi_k(t)\,dt = 1
\]
for any $k$. On the other hand, by \eqref{rearra}
\begin{align*}
\|\phi_k\|^2_{W^{2,2}(M)} 
&\geq \omega_n\int_{-\infty}^{+\infty}
\left[j(t)^{n-1}+ (n-1) j(t)^{n-3}(j'(t))^2\right](\p_t \hat \phi_k(t))^2 \,dt.
\end{align*}
Suppose now that 
\begin{equation}\label{jderiv}
\int_{-\infty}^{+\infty}
\left[j(t)^{n-1}+ (n-1) j(t)^{n-3}(j'(t))^2\right]^{-1} \,dt =: \mathfrak J <+\infty.
\end{equation}
Then, by Cauchy-Schwarz inequality,
\begin{align}\label{CS}
\|\phi_k\|^2_{W^{2,2}(M)} 
&\geq \omega_n \frac{\left(\int_{-\infty}^{\infty} \p_t\hat\phi_k(t)\,dt\right)^2}{\int_{-\infty}^{+\infty}
	\left[j(t)^{n-1}+ (n-1) j(t)^{n-3}(j'(t))^2\right]^{-1} \,dt}\\
&=\omega_n\mathfrak J^{-1}>0\nonumber
\end{align}
independently of $k$.
Accordingly, to conclude the proof, we have to construct a warping function $j:\R\to (0,\infty)$ which has fast enough decay in order to satisfy \eqref{jvol}, but which is sufficiently oscillating  in order to satisfy also \eqref{jderiv}.
 
Note that it is enough to define $j$ on $[1,+\infty)$ such that $j\in L^1([1,+\infty))$ and  
\begin{equation}\label{jderivplus}
\int_{1}^{+\infty}
\left[j(t)^{n-1}+ (n-1) j(t)^{n-3}(j'(t))^2\right]^{-1} \,dt =: \mathfrak J_+ <+\infty.
\end{equation}
The definition of $j$ will be then extended to $(-\infty,-1]\cup[1,+\infty)$ by symmetry, i.e. setting $j(t)=j(-t)$, and finally to the whole $\R$ taking any positive smooth extension in $[-1,1]$.

Choose the warping function $j$ in such a way that for all $t>1$ it holds \[
j(t)=\psi(t)t^{-\frac 1{n-1}-\epsilon},
\]
where $\psi\in C^\infty([1,\infty))$ is a function, to be chosen later, which satisfies 
\begin{equation}\label{eq psi}
1\leq\psi\leq 2.
\end{equation}
In particular
\[
j^{n-1}(t)\leq 2^{n-1}|t|^{-1-\epsilon(n-1)}
\]
for $|t|>1$, so that $j^{n-1}
\in L^1(\R)$,
and \eqref{jvol} is satisfied. 

In order to construct $\psi$, we define first a continuous piecewise-smooth function $\tilde\psi:[1,+\infty)\to[0,+\infty)$ by 
\[
\tilde\psi(t):=\begin{cases} 
t^{2^+}-\lfloor t^{2^+} \rfloor+1,&\text{if} \lfloor t^{2^+} \rfloor\text{ is even},\\
-t^{2^+}+\lfloor t^{2^+} \rfloor +2,&\text{if} \lfloor t^{2^+} \rfloor \text{ is odd}
\end{cases}
\]
where we have set $2^+:= 2+ \epsilon n/2$, and $\epsilon>0$ is a fixed small constant. Define $\tilde j(t):=\tilde\psi(t)t^{-\frac 1{n-1}-\epsilon}$. Outside the (null-measure) singular set 
\[
\mathcal S:= \{t\in[1,+\infty)\,:\,t^{2^+}=m,\ m\in\mathbb N\}=\{t_m\}_{m=1}^\infty,\]
with $t_m:=m^{1/{2^+}}$, we can compute
\begin{align*}
\left|\tilde j'
\right|&=\left|\frac{\tilde\psi'(t)t-\left(\frac{1}{n-1}+\epsilon\right)\tilde\psi(t)}{t^{\frac{n}{n-1}+\epsilon}}\right|\\
&\ge t^{-\frac{n}{n-1}-\epsilon} \left\{\left|\tilde\psi'(t)t\right|-\left(\frac{1}{n-1}+\epsilon\right)|\tilde\psi(t)|\right\}\\
&\ge t^{-\frac{n}{n-1}-\epsilon} \left\{\left|\tilde\psi'(t)t\right|-2\left(\frac{1}{n-1}+\epsilon\right)\right\}\\
&\ge 2^+ t^{2^+-\frac{n}{n-1}-\epsilon}-2\left(\frac{1}{n-1}+\epsilon\right)t^{-\frac{n}{n-1}-\epsilon}
\\
&\ge t^{\frac{n-2}{n-1}+\epsilon(\frac{n-2}{2})} \end{align*}
for $t\geq t_0:= \left(2/(n-1)+2\epsilon\right)^{1/2^+}$.
In particular 
\begin{align}\label{frakjplus}
\tilde {\mathfrak J}_+:&=\int_{1}^{+\infty}
\left[\tilde j(t)^{n-1}+ (n-1) \tilde j(t)^{n-3}(\tilde j'(t))^2\right]^{-1} \,dt\\
&\leq \frac 1{n-1}\int_{1}^{+\infty}
\left[ \tilde j(t)^{n-3}(\tilde j'(t))^2\right]^{-1} \,dt\nonumber\\
&\leq \frac 1{n-1}\int_{1}^{+\infty}
\left[ t^{-\frac {n-3}{n-1}-\epsilon(n-3)}t^{\frac{2n-4}{n-1}+\epsilon(n-2)}\right]^{-1} \,dt\nonumber\\
&\leq \frac 1{n-1}\int_{1}^{+\infty}
t^{-1-\epsilon} \,dt\nonumber\\
& <+\infty.\nonumber
\end{align}
Since we want a smooth Riemannian metric, to conclude the proof we have to smooth away $\tilde\psi$ (hence $\tilde j$) in a neighbourhood of the singular set $\mathcal S$ in such a way that \eqref{jderivplus} and \eqref{eq psi} are preserved. To this end, it is enough to take a smooth uniform approximations of $\psi$ in a sufficiently small neighbourhood of $\mathcal S$. 

 For any integer $m\geq 1$, set $\eta_m:=m^{-(3+n\epsilon)/2}$.
Define 
\[
A:=[1,\infty)\cap\left(\cup_{m\geq 1}[t_m-\eta_m,t_m+\eta_m]\right),\qquad B:=[1,\infty)\setminus A.\]
 Set $\psi(t):=\tilde \psi(t)$ for $t\in B$ and extend $\psi$ to the whole $[1,\infty)$ by taking any $C^\infty$ prolongation of $\psi|_B$ which satisfies \eqref{eq psi}. We write 
\begin{align*}
{\mathfrak J}_+:=&\int_{1}^{+\infty}
\left[ j(t)^{n-1}+ (n-1) j(t)^{n-3}( j'(t))^2\right]^{-1} \,dt\\
&=\int_A
\left[ j(t)^{n-1}+ (n-1) j(t)^{n-3}( j'(t))^2\right]^{-1} \,dt
+\int_B
\left[\tilde j(t)^{n-1}+ (n-1) \tilde j(t)^{n-3}(\tilde j'(t))^2\right]^{-1} \,dt\\
 &\leq \int_A
j(t)^{1-n} \,dt
+\tilde {\mathfrak J}_+\\
&\leq \int_A
t^{1+\epsilon(n-1)}\,dt
+\tilde {\mathfrak J}_+
\\
 &\leq \sum_{m=1}^\infty \int_{t_m-\eta_m}^{t_m+\eta_m}
t^{1+\epsilon(n-1)}\,dt
+\tilde {\mathfrak J}_+\\
 &\leq \sum_{m=1}^\infty 2\eta_m
(t_m+\eta_m)^{1+\epsilon(n-1)}
+\tilde {\mathfrak J}_+.
\end{align*}

Since $\eta_m=o(t_m)$, we have 
\begin{align*}
2\eta_m
(t_m+\eta_m)^{1+\epsilon(n-1)}
&\leq 2\eta_m (2t_m)^{1+\epsilon(n-1)}\\ &
\leq 2m^{-\frac{3+\epsilon n}{2}} 2^{1+\epsilon(n-1)}m^{\frac{2+2\epsilon(n-1)}{4+\epsilon n}} 	\\
&=2^{2+\epsilon(n-1)}m^{\frac{4+4\epsilon(n-1)-12-4\epsilon n-3\epsilon n- \epsilon^2 n^2}{8+2\epsilon n}}\\
&=2^{2+\epsilon(n-1)}m^{\frac{-8-4\epsilon -3\epsilon n- \epsilon^2 n^2}{8+2\epsilon n}}\\
&\leq 2^{2+\epsilon(n-1)}m^{-\frac{8+3\epsilon n}{8+2\epsilon n}}.
\end{align*}
Then $\sum_{m=1}^\infty 2\eta_m
(t_m+\eta_m)^{1+\epsilon(n-1)}<\infty$, and together with \eqref{frakjplus}, this gives $\mathfrak J_+<\infty$ and conclude the proof for $p=2$. 


\begin{remark}\label{rmk_CZ}{\rm
We finish this section observing that the manifold $(M,g)$ constructed in the proof of Theorem \ref{th_main} gives also a counterexample to the validity of global Calder\'on-Zygmund inequalities on Riemannian manifolds (see \cite{GuneysuPigola} and \cite{Li} for previous examples). I'm grateful to S. Pigola for pointing out this fact to me.\\
Indeed, suppose by contradiction that the Calder\'on-Zygmund inequality 
\begin{equation}\label{CZ}
\| \Hess_g u \|_{L^p}\leq C \{\|\Delta_g u\|_{L^p} + \|u\|_{L^p} \},\qquad \forall\,u\in C^\infty_c(M)
\end{equation}
holds on $M$ for some $p\geq 2$. Let $f\in H_0^{2,p}(M)$ and $\{f_k\}\subset C^\infty_c$ such that $f_k\to f$ in $H^{2,p}(M)$. Then $\{f_k\}$ is a Cauchy sequence in $H^{2,p}(M)$, and thus it is also Cauchy in $W^{2,p}(M)$ because of \eqref{CZ} and \cite[Theorem 2]{GuneysuPigola_AMPA}. Since $W^{2,p}(M)$ is a Banach space, $\{f_k\}$ converges some $\bar f\in W^{2,p}_0(M)$. But $\bar f = f$ because of the continuous embedding $W^{2,p}(M)\subseteq H^{2,p}(M)$. Accordingly 	\[H^{2,p}(M)=H_0^{2,p}(M)\subseteq W^{2,p}_0(M)\subseteq W^{2,p}(M)\subseteq H^{2,p}(M),\] which contradicts $W^{2,p}_0(M)\subsetneq W^{2,p}(M)$.  
}\end{remark}

\section{An example of $W^{2,p}(M)\subsetneq H^{2,p}(M)$}\label{app}
In this final section, for any $p\in [1,\infty)$ and any dimension $n\geq 2$, we give an example of complete non-compact $n$-dimensional Riemannian manifolds $(M,g)$ on which $H^{2,p}(M)\setminus W^{2,p}(M)$ is non-empty. The construction is essentially due to B. G\"uneysu and S. Pigola (for $n=p=2$), \cite[Theorem 4.6]{GuneysuPigola}, and to S. Li (for $n\geq 3$ and $p\neq 2$), \cite{Li}, and was originally used to give an example of a manifold which does not support an $L^p$-Calder\'on-Zygmund inequalities. 

\begin{proposition}
	For any $p\in [1,\infty)$ and any dimension $n\geq 2$, there exists a complete non-compact $n$-dimensional Riemannian manifold $(M,g)$ such that the set $H^{2,p}(M) \setminus W^{2,p}(M) $ is non-empty.
	\end{proposition}

As pointed out in the introduction, it is worth to notice that this result is possibly not new. At least in the analogous case of Sobolev spaces of $1$-forms a counterexample due to J. Dodziuk can be found in \cite[p. 67]{dodziuk}.

\begin{proof}
	Fix $p\in [1,\infty)$ and $n\geq 2$. In \cite{Li}, the author constructed a Riemannian manifold $(M,g)$ and a family of functions $\{u_k\}_{k\in\mathbb N}\subset C^\infty_c(M)$ such that \begin{itemize}
		\item the functions $u_k$'s are supported on intervals with pairwise disjoint interiors;
		\item $\|u_k\|_{L^p(M)}^p\leq C_{n,p} e^{2k}$;
	\item $\|\Delta u_k\|_{L^p(M)}^p\leq C_{n,p} e^{-2(p-1)k}$;
\item $\|\Hess\, u_k\|_{L^p(M)}^p\geq C_{n,p}^{-1} e^{k^{1000}}$ for any sufficiently large $k$.
\end{itemize} Here $C_{n,p}$ is a positive constant depending only on $n$ and $p$. 
Define $F\dot = \sum_{k=0}^\infty e^{-3k}u_k$. The sum converges, since it is locally finite. Moreover
\[
\|F\|_{L^p(M)}^p\leq C_{n,p} \sum_{k=0}^\infty e^{-3k}e^{2k}= C_{n,p} \frac e{e-1}
\]
and 
\[
\|\Delta F\|_{L^p(M)}^p\leq C_{n,p} \sum_{k=0}^\infty e^{-3k}e^{-2(p-1)k}\leq C_{n,p} \frac e{e-1},
\]
so that $F\in H^{2,p}(M)$. However
\[
\|\Hess\, F\|_{L^p(M)}^p\geq C_{n,p}^{-1} \sum_{k=0}^\infty e^{-3k}e^{k^{1000}}=+\infty,
\]
so that $F\not\in W^{2,p}(M)$.

	\end{proof} 

\begin{acknowledgement*}
	I would like to thank S. Pigola for useful comments on the first version of this paper. In particular I'm grateful to him for pointing me out Remark \ref{rmk_CZ}.\\ The author is member of the "Gruppo Nazionale per l'Analisi
	Matematica, la Probabilità e le loro Applicazioni" (GNAMPA) of the Istituto
	Nazionale di Alta Matematica (INdAM). 
\end{acknowledgement*}

\bibliographystyle{amsplain}
\bibliography{CutOffs}

  \makeatletter
    \providecommand\@dotsep{5}
  \makeatother

\end{document}